\documentclass[11pt,a4paper]{article}
\usepackage{ifthen,latexsym,amssymb,amsmath,fixmath,bbm}

\newcommand{\C}[1]{{\protect\cal #1}}
\newcommand{\B}[1]{\mathbold{#1}}
\newcommand{\I}[1]{{\mathbb #1}}

\newcommand{\ceil}[1]{\lceil #1\rceil}

\newcommand{\floor}[1]{\lfloor #1\rfloor}
\newcommand{\me}{{\mathrm e}}
\renewcommand{\mid}{:}
\renewcommand{\ge}{\geqslant}
\renewcommand{\le}{\leqslant}
\renewcommand{\preceq}{\preccurlyeq}
\renewcommand{\succeq}{\succcurlyeq}

\newif\ifnotesw\noteswtrue

\newcommand{\hide}[1]{}

\newcommand{\beq}[1]{\begin{equation}\label{#1}}
\newcommand{\eeq}{\end{equation}}

\newtheorem{theorem}{Theorem}
\newtheorem{lemma}[theorem]{Lemma}

\newtheorem{corollary}[theorem]{Corollary}

\newcommand{\bpf}[1][Proof.]{\smallskip\noindent{\it #1} }
\newcommand{\qed}{\nolinebreak\mbox{\hspace{5 true pt}%
  \rule[-0.85 true pt]{3.9 true pt}{8.1 true pt}}}
\newcommand{\epf}{\qed \medskip}

\newcommand{\OPdata}{Oleg Pikhurko\footnote{Supported by ERC
grant~306493 and EPSRC grant~EP/K012045/1.
 }\\
Mathematics Institute and DIMAP\\
University of Warwick\\
Coventry CV4 7AL, UK}

\begin{document}


\newcommand{\IS}{\I S} 
\newcommand{\ex}{\mathrm{ex}}
\newcommand{\FinPi}[1]{\Pi^{(#1)}_{\mathrm{fin}}}
\newcommand{\InfPi}[1]{\Pi_{\infty}^{(#1)}}
\newcommand{\branch}[2]{\mathrm{br}_{#1}(#2)}
\newcommand{\blow}[2]{#1(\!(#2)\!)}
\newcommand{\density}[1]{\Lambda_{#1}}
\newcommand{\multiset}[1]{\{\hspace{-0.25em}\{\hspace{0.1em}#1\hspace{0.1em}\}\hspace{-0.25em}\}}
\newcommand{\E}[1]{\mathbf{E}\left[\hspace{0.5mm}#1\hspace{0.5mm}\right]}
\newcommand{\ind}{{}_{\mathrm{ind}}}
\newcommand{\LIM}[1]{\mathrm{LIM}^{(#1)}}
\newcommand{\CG}{\C G}
\newcommand{\cl}[1]{\mathrm{cl}(#1)}
\newcommand{\rep}[2]{{#1}^{(#2)}}
\newcommand{\order}[1]{\langle #1\rangle}
\newcommand{\ins}[2]{\order{\B{#1},#2}}
\newcommand{\minus}[2]{#1\setminus #2}

\author{\OPdata}
\title{The maximal length of a gap between\\ $r$-graph Tur\'an densities}
\maketitle

\begin{abstract}
The \emph{Tur\'an density} $\pi(\C F)$ of a family $\C F$ of
$r$-graphs is the limit as $n\to\infty$ 
of the maximum edge density 
of an $\C F$-free $r$-graph on $n$ vertices.  
Erd\H os [%
Israel J. Math \textbf{2} (1964) 183--190] proved
that no Tur\'an density can lie in the open interval $(0,r!/r^r)$. Here we show that any other
open subinterval of $[0,1]$ avoiding Tur\'an densities has strictly smaller 
length. In particular, this implies a conjecture of Grosu~[%
E-print \texttt{arXiv:1403.4653v1}, 2014].
\end{abstract}

\section{Introduction}

Let $\C F$ be a (possibly infinite) family of \emph{$r$-graphs} (that is,
$r$-uniform set systems). We call elements of $\C F$ \emph{forbidden}. An
$r$-graph $G$ is \emph{$\C F$-free} if no member $F\in \C F$ is 
a subgraph of $G$, that is, we cannot obtain $F$
by deleting some vertices and edges from $G$. The \emph{Tur\'an
function $\ex(n,\C F)$} is the maximum number of edges that an $\C F$-free
$r$-graph on $n$ vertices can have. This is one of the central questions of
extremal combinatorics that goes back to the fundamental
paper of Tur\'an~\cite{turan:41}. We refer the reader to the
surveys of the Tur\'an function by
F\"uredi~\cite{furedi:91}, Keevash~\cite{keevash:11},
and Sidorenko~\cite{sidorenko:95}.

As it was observed by Katona, Nemetz, and
Simonovits~\cite{katona+nemetz+simonovits:64}, 
the limit
 $$
 \pi(\C F):=\lim_{n\to\infty} \frac{\ex(n,\C F)}{{n\choose k}}
 $$
 exists. It is called the \emph{Tur\'an density} of $\C F$. Let $\InfPi{r}$
consist of all possible Tur\'an densities of $r$-graph families and let
 $\FinPi{r}$ be the set of all possible Tur\'an densities when \textit{finitely many}
$r$-graphs are forbidden.  It is convenient to allow empty forbidden families, so that 
$1$ is also a Tur\'an density. Clearly, $\FinPi{r}\subseteq \InfPi{r}$. 
A result of Brown and Simonovits~\cite[Theorem~1]{brown+simonovits:84} implies that 
the topological closure $\cl{\FinPi{r}}$ of $\FinPi r$ contains $\InfPi r$ while the converse inclusion
was established in~\cite[Proposition~1]{pikhurko:14}; thus
 \begin{equation}\label{eq:cl}
 \InfPi r=\cl{\FinPi{r}},\quad \mbox{for every integer $r\ge 2$}.
 \end{equation}

For $r=2$, the celebrated Erd\H os-Stone-Simonovits 
Theorem~\cite{erdos+simonovits:66,erdos+stone:46} determines 
the Tur\'an density for every family $\C F$. In particular, we have
 \begin{equation}\label{eq:PI2}
 \FinPi{2}=\InfPi{2}=\left\{\frac{m-1}m\mid m=1,2,3,\dots,\infty\right\}.
 \end{equation}

Unfortunately, the Tur\'an function for
\emph{hypergraphs} (that is, $r$-graphs with 
$r\ge 3$) is much more difficult to analyse and many problems
(even rather basic ones) are wide open. 

Fix some $r\ge 2$. A \emph{gap}
is an open interval $(a,b)\subseteq (0,1)$ that is disjoint
from $\InfPi{r}$ (which, by~\eqref{eq:cl}, is equivalent to being
disjoint from $\FinPi{r}$). Here we consider $g_r$, the maximal possible length of a gap.
In other words, $g_r$ is
the maximal $g$ such that there is a real $a$ with $(a,a+g)\subseteq (0,1)\setminus
\InfPi{r}$.  For example, \eqref{eq:PI2} implies that $g_2=1/2$. Erd\H os~\cite{erdos:64} proved
that $(0,r!/r^r)$ is a gap; in particular, $g_r\ge r!/r^r$. Here we show that 
this is equality and every other gap has strictly smaller length.

\begin{theorem}\label{th:gap} For every $r\ge 3$, we have that $g_r=r!/r^r$ and, furthermore, $(0,r!/r^r)$ is the only gap of length $r!/r^r$ for $r$-graphs.
\end{theorem}

In particular we obtain the following result that was conjectured by Grosu~\cite[Conjecture 10]{grosu:densities}.

\begin{corollary}\label{cr:grosu} The union of $r$-graph Tur\'an densities over all
$r\ge 2$ is dense in $[0,1]$, that is, $\cl{\cup_{r=2}^\infty \InfPi r}=[0,1]$.\qed
 \end{corollary}

The question whether the set $\InfPi{r}$  is
a well-ordered subset of $([0,1],<)$ for $r\ge 3$ was a famous \$1000 problem of Erd\H os that was answered in the negative by Frankl and R\"odl~\cite{frankl+rodl:84}. Despite a number of results
that followed~\cite{frankl+rodl:84}, very little is known about other gaps in $\InfPi{r}$ for $r\ge3$. For example, let $g_r'$
be the the second largest gap length, that is, the maximum $g\ge 0$ such that $(a,a+g)\subseteq (r!/r^r,1)\setminus\InfPi r$
for some $a$. The computer-generated proof of 
Baber and Talbot~\cite{baber+talbot:11} implies that $g_3'\ge 0.0017$. 
However, not for a single $r\ge 4$ is it known, for example, whether $g_r'$ is zero (i.e.\ whether
$\InfPi{r}$ is dense in $[r!/r^r,1]$).  

This paper is organised as follows. In Section~\ref{prelim} we give some definitions and auxiliary results. 
Theorem~\ref{th:gap} is proved in Section~\ref{first}. We give another proof of Corollary~\ref{cr:grosu} in Section~\ref{second}. Although the latter proof is not strong enough
to prove Theorem~\ref{th:gap}, its advantage is that it produces explicit elements of $\FinPi r$ (as opposed to 
the implicit values of certain maximisation problems returned by the proof in Section~\ref{first}). 
So we include both proofs here, even though the second one is longer. 

\section{Preliminaries}\label{prelim}

For $n\in \I N$, define $[n]:=\{1,\dots,n\}$. For reals $a\le b$, let $(a,b)$ and $[a,b]$ be respectively open
and closed intervals of reals with endpoints $a$ and $b$. The \emph{standard $(m-1)$-dimensional simplex} is
 $$
 \IS_m:=\{\B x\in \I R^m\mid x_1+\dots+x_m=1,\ \forall\, i\in[m]\ x_i\ge 0\}.
 $$

An \emph{$r$-pattern} is a
collection $P$ of $r$-multisets on $[m]$, for some $m\in\I N$.
(By an \emph{$r$-multiset} we mean an unordered collection
of $r$ elements with repetitions allowed.) Let $V_1,\dots,V_m$ 
be disjoint sets and let $V=V_1\cup\dots\cup V_m$. The \emph{profile}
of an $r$-set $X\subseteq V$ (with respect to $V_1,\dots,V_m$) is
the $r$-multiset on $[m]$ that contains $i\in [m]$
with multiplicity $|X\cap V_i|$. For an $r$-multiset $Y$ on~$[m]$, let 
$\blow{Y}{V_1,\dots,V_m}$ consist of all $r$-subsets of $V$ whose profile is
$Y$. We call this $r$-graph the \emph{blow-up of $Y$} (with respect to $V_1,\dots,V_m$) and
the $r$-graph
 $$
 \blow{P}{V_1,\dots,V_m}:=\bigcup_{Y\in P} \blow{Y}{V_1,\dots,V_m}
 $$ 
 is called the \emph{blow-up of $P$}. 
Let the \emph{Lagrange polynomial} of $P$ be 
 $$
 \lambda_P(x_1,\dots,x_m):=r!\,\sum_{D\in P}\; \prod_{i=1}^m\;
\frac{x_i^{D(i)}}{D(i)!}\in\I R[x_1,\dots,x_m], 
 $$
 where $D(i)$ denotes the multiplicity of $i$ in $D$.
This definition is motivated by the fact that, 
for every partition $[n]=V_1\cup\dots\cup V_m$, we have that
 $$
|\blow{P}{V_1,\dots,V_m}|
= \lambda_P\left(\frac{|V_1|}n,\dots,\frac{|V_m|}n\right)\times {n\choose r}
+O(n^{r-1}),\qquad  \mbox{as $n\to\infty$}.
 $$

For example, if $r=3$, $m=3$, and $P$ consists of multisets $\{1,1,2\}$ and $\{1,2,3\}$,
then $\blow{P}{V_1,\dots,V_m}$ contains all triples that have two vertices in $V_1$ and one vertex in $V_2$ plus all triples with exactly one vertex in each part; here $\lambda_P(x_1,x_2,x_3)=3 x_1^2x_2+6x_1x_2x_3$.

 Let the \emph{Lagrangian of $P$} be $\Lambda_P:=\max\{\lambda_P(\B x)\mid \B x\in \IS_m\}$,  the maximum value of the polynomial 
$\lambda_P$ on the compact set~$\IS_m$. One obvious connection of this parameter to $r$-graph Tur\'an
densities is that, if each blow-up of $P$ is $\C F$-free, then $\pi(\C F)\ge \Lambda_P$.  Also, it is not hard
to show that $\Lambda_P=\pi(\C F)$, where $\C F$ consists
of all $r$-graphs $F$ such that every blow-up of $P$ is $F$-free; thus $\Lambda_P\in \InfPi r$.
As shown in~\cite[Theorem~3]{pikhurko:14}, we have in fact that
  \beq{eq:Lambda} 
  \Lambda_P\in \FinPi r,\quad\mbox{for every $r$-pattern $P$}.
  \eeq


We will use
a special case of Muirhead's inequality (see e.g.~\cite[Theorem 45]{hardy+littlewood+polya:i}) which states that, for any $0\le i<j\le k$, we have
 \begin{equation}\label{eq:Muirhead}
 x^{k+i}y^{k-i}+x^{k-i}y^{k+i}\le x^{k+j}y^{k-j}+x^{k-j}y^{k+j},\quad \mbox{for }x,y\ge 0.
 \end{equation}

\section{Proof of Theorem~\ref{th:gap}}\label{first}

Let $r\ge 3$. Fix a sufficiently large integer $m=m(r)$
so that $r!{m\choose r}/m^{r}> 1-r!/r^r$. Consider $r$-graphs $G_0,\dots,G_{{m\choose r}}$ on $[m]$ such that $G_0$ has no edges and, for $i=1,\dots,{m\choose r}$, the $r$-graph $G_i$ is obtained
from $G_{i-1}$ by adding a new edge. In other words, we enumerate all $r$-subsets of $[m]$ as $R_1,\dots,R_{{m\choose r}}$
and let $G_i:=\{R_1,\dots,R_i\}$. Let 
 $$
\lambda_i(\B x):=\lambda_{G_i}(\B x)=r!\sum_{D\in G_i} \prod_{j\in D} x_j,
 $$ 
 be the Lagrange polynomial of $G_i$ and $\Lambda_i:=\Lambda_{G_i}$ be its Lagrangian, where we view
$G_i$ as an $r$-pattern. Since $G_{i-1}\subseteq G_i$, we have that
$\Lambda_{i-1}\le \Lambda_i$. 

We claim that for every $i\in [\,{m\choose r}\,]$
 \beq{eq:aim3}
 \Lambda_i-\Lambda_{i-1}\le r!/r^r.
 \eeq
 Indeed, pick $\B x\in\IS_m$ with $\Lambda_i=\lambda_i(\B x)$. Let $R_i=\{u_1,\dots,u_r\}$. When we remove the term $r!\, x_{u_1}\ldots x_{u_r}$
from $\lambda_i(\B x)$, we get the evaluation of $\lambda_{i-1}$ on $\B x\in\IS_m$. By definition, $\Lambda_{i-1}\ge \lambda_{i-1}(\B x)$. 
Also, since $x_{u_1}+\ldots+x_{u_r}\le 1$, we have $x_{u_1}\ldots x_{u_r}\le r^{-r}$ by the Geometric-Arithmetic Mean Inequality. Thus we obtain the stated bound:
 $$
 \Lambda_i=\lambda_i(\B x)=\lambda_{i-1}(\B x)+r!\,x_{u_1}\ldots x_{u_r}\le \Lambda_{i-1}+r!/r^r.
 $$ 

Also, we have $\Lambda_{{m\choose r}}\ge \lambda_{{m\choose r}}(\frac1m,\dots,\frac1m)= r!{m\choose r}/m^{r}> 1-r!/r^r$. This and~\eqref{eq:Lambda} imply that $g_r\le r!/r^r$, while
the result of Erd\H os~\cite{erdos:64} gives the converse inequality. Also, if we 
have equality in \eqref{eq:aim3}, then necessarily $x_{u_1}=\dots=x_{u_r}=1/r$ and thus $\Lambda_{i-1}=0$,
implying the uniqueness part of  Theorem~\ref{th:gap}.

\section{Alternative proof of Corollary~\ref{cr:grosu}}\label{second}

For integers $r,s\ge 2$, let $\C P_{r,s}$ consist of ordered $s$-tuples $(r_1,\dots,r_s)$ of non-negative integers such that $r_1\ge\dots\ge r_s$ and $r_1+\dots+r_s=r$. 
This set admits a partial order, where $\B x\succeq \B y$ if $\sum_{i=1}^k x_i\ge \sum_{i=1}^k y_i$ 
for every $k\in [s]$.
For example, the (unique) maximal element is $(r,0,\dots,0)$ and the (unique) minimal element is $(\ceil{r/s},\dots,\floor{r/s})$. 

Let $A\subseteq
\C P_{r,s}$. The set $A$ is called \emph{down-closed} if $\B y\in A$ whenever $\B x\in A$ and $\B x\succeq \B y$. 
Let $G_A$ consist of all $r$-multisets $X$ on $[s]$ such that the mutiplicities of $X$ satisfy $\order{X(1),\dots,X(s)}\in A$,
where $\order{\B x}$ denotes the non-increasing ordering of a vector~$\B x$. Also, we use shortcuts $\lambda_{A}:=\lambda_{G_A}$ and $\Lambda_A:=\Lambda_{G_A}$.

\begin{lemma}\label{lm:opt} Let $r,s\ge 2$. If $A\subseteq \C P_{r,s}$ is down-closed, then
$\Lambda_A=\lambda_A(\frac 1s,\dots,\frac 1s)$.\end{lemma}

\bpf We use induction on $s$. 

First, we prove the base case $s=2$. Let $k:=r/2$. For $h\ge 0$, let $I_h$ consist of all integer translates of $k$ whose absolute value is at most $h$, that is, 
$I_h:= (\I Z+k)\cap [-h,h]$. Also, let $I_h^+:=I_h\cap [0,h]$. (These definitions will allow us to deal with the cases of even and odd
$r$ uniformly.) For example, $\C P_{r,2}=\{(k+ i,k- i)\mid i\in I_k^+\}$. 

Take a down-closed set $A\subseteq \C P_{r,2}$. It consists of pairs $(k+i,k-i)$ with $i\in I_h^+$ for some $h$. 
By the homogeneity of the polynomials involved, the required inequality can be rewritten as 
 \beq{eq:aim}
 \sum_{i\in I_h} {2k\choose k+i}\left(\frac{x+y}2\right)^{2k}-\sum_{i\in I_h} {2k\choose k+i}x^{k+i}y^{k-i}\ge 0,\quad \mbox{for }x,y\ge 0.
 \eeq

 We will apply the so-called \emph{bunching method} where we try to write the desired
inequality as a positive linear combination of Muirhead's inequalities~\eqref{eq:Muirhead}. If 
$j\in I_h$, then the coefficient in front of $x^{k+j}y^{k-j}$ in (\ref{eq:aim})
is
 $$
 2^{-2k}{2k\choose k+j}\sum_{i\in I_h} {2k\choose k+i}- {2k\choose k+j}\le 0.
 $$
 If $j\in I_k\setminus I_h$, then the coefficient is
 $
 2^{-2k}{2k\choose k+j}\sum_{i\in I_h} {2k\choose k+i}\ge 0.
 $
 Thus, if we group (\ref{eq:aim}) into terms $x^{k+j}y^{k-j}+x^{k-j}y^{k+j}$, 
then we get non-positive coefficients for $0\le j\le h$ followed by non-negative
coefficients for $j>h$. Also, the total sum of coefficients is zero because
(\ref{eq:aim}) becomes equality for $x=y=1$. Thus we can ``bunch'' 
$I_h$-terms with $(I_k\setminus I_h)$-terms and use~\eqref{eq:Muirhead} to derive the desired inequality (\ref{eq:aim}).
This proves the case $s=2$.

Now, let $s\ge 3$ and suppose that we have proved the lemma for $s-1$ (and all $r$). The function
$\lambda_A$ is a continuous function on the compact set $\IS_s$. Let it attain its maximum
on some $\B x\in\IS_s$. If there is more than one choice, then choose $\B x$ so that
$\Delta:=\sum_{i\not=j} |x_i-x_j|$ is minimised.  Suppose that $\Delta\not=0$, say $x_1\not=x_2$.
Note that $\lambda_A$ is a homogeneous polynomial of degree $r$, 
and the coefficient at $x_1^{r_1}\ldots x_s^{r_s}$ is ${r\choose r_1,\dots,r_s}$ if the ordering $\order{\B r}$ of $\B r$ is in $A$ and 0 otherwise.

Fix $j\in \{0,\dots,r\}$. If we collect all terms in front of $x_s^j$, we get
 $$
 \sum_{\ins{r}{j}\in A\atop r_1+\dots+r_{s-1}=r-j} {r\choose r_1,\dots,r_{s-1},j} \prod_{i=1}^{s-1} x_i^{r_i}= 
{r\choose j}\, \lambda_{\minus{A}{j}}(x_1,\dots,x_{s-1}),
 $$
 where $\ins{y}{j}$ is 
obtained from $\B y$  by inserting $j$ and ordering the obtained sequence, while $\minus{A}{j}$ consists of those $\B y\in \C P_{r-j,s-1}$ such that $\ins{y}{j}\in A$.

Let us show that $\minus{A}{j}\subseteq \C P_{r-j,s-1}$ is down-closed. Take arbitrary $\B z\in \minus{A}{j}$
and $\B y\preceq \B z$. We have to show that $\B y\in \minus{A}{j}$. 
Since $A\ni \ins{z}{j}$ is down-closed, it is enough to show that $\ins{z}{j}\succeq \ins{y}{j}$. We have to compare the sums of the first $i$ terms of $\ins{z}{j}$
and of $\ins{y}{j}$. A problem could arise only if the new entry $j$ was included  into these terms for $\ins{y}{j}$, say as the term number $h\le i$, but not for $\ins{z}{j}$. Since $\B z\succeq \B y$, we have that $\sum_{f=1}^{h-1} z_f\ge \sum_{f=1}^{h-1} y_f$ (and these are also the initial sums for $\ins{z}{j}$ and $\ins{y}{j}$).
Furthermore, each of the subsequent $i-(h-1)$ entries is at least $j$ for $\ins{z}{j}$ and at most $j$ for $\ins{y}{j}$.
It follows that $\ins{z}{j}\succeq \ins{y}{j}$. Thus $\minus{A}{j}$ is down-closed, as claimed.

By the induction assumption (and since $\lambda_{\minus{A}{j}}$ is a homogeneous polynomial), we have that $\lambda_{\minus{A}{j}}(x_1,\dots,x_{s-1})\le \lambda_{\minus{A}{j}}(\frac{1-x_s}{s-1},\dots,\frac{1-x_s}{s-1})$. Thus
 $$
 \Lambda_A=\lambda_A(\B x)=\sum_{j=0}^r {r\choose j}\, \lambda_{\minus{A}{j}}(x_1,\dots,x_{s-1})\, x_s^j \le \lambda_{A}\left(\frac{1-x_s}{s-1},\dots,\frac{1-x_s}{s-1},x_s\right).
 $$
 Clearly, the sum $\sum_{i=1}^{s-1}|x_s-x_i|$ does not increase if we replace each of $x_1,\dots,x_{s-1}$ by their
arithmetic mean $(1-x_{s})/(s-1)$. Since $x_1\not=x_2$, we have found another optimal element of
$\IS_s$ with strictly smaller $\Delta$, a contradiction. The lemma is proved.\epf

Fix some enumeration $\C P_{r,r}=\{R_1,\dots,R_m\}$ such that if $R_i\succeq R_j$ then $i\ge j$. For $j\in \{0,\dots,m\}$, let $A_j:=\{\,R_i\mid i\in [j]\,\}$.
Thus, for example, $A_0=\emptyset$ and $A_m=\C P_{r,r}$. By~\eqref{eq:Lambda}, $\FinPi{r}$ contains
all of the following numbers:
 $$
 0=\Lambda_{A_0}\le \Lambda_{A_1}\le \dots\le \Lambda_{A_m}=1.
 $$

\hide{
Unfortunately, this sequence does have gaps strictly larger than $r!/r^r$. For example, $\Lambda_{A_1}-\Lambda_{A_0}= {r\choose 2}r(r-1)(r-2)! r^{-r}={r\choose 2}\, \frac{r!}{r^r}$: using the balls-urns interpretation, we first choose $2$ balls that go into the same urn, then the two special urns (for 2 and for 0 balls) and then place the other $r-2$ balls bijectively into
$r-2$ remaining urns.
}%

Let us
show that $\max\{\,\Lambda_{A_i}-\Lambda_{A_{i-1}}\mid i\in [m]\,\}=o(1)$ 
as $r\to\infty$. By definition, each $A_j\subseteq \C P_{r,r}$ is down-closed. Thus, by Lemma~\ref{lm:opt} the difference 
 $
 \Lambda_{A_i}-\Lambda_{A_{i-1}}
 $
is
the probability that, when $r$ balls are uniformly and independently distributed into $r$ urns, the
ordered ball distribution is given by $R_i$. 
Expose the first $r-m$ balls, where, for example, $m:=\lfloor \log r\rfloor$. 
Let $k$ be the number of empty cells. Its expected value is $r (1-1/r)^{r-m}=(\me^{-1}+o(1))\,r$.
By Azuma's inequality (see e.g.\ \cite[Theorem~7.2.1]{alon+spencer:pm}), we have
\emph{whp} (i.e.\ with probability $1-o(1)$ as $r\to\infty$) that $k$ is in $I:=[r/4,3r/4]$. Assume that $k\in I$ and expose the
remaining $m$ balls. Let $J$ be the number of balls that land inside the $k$  cells that were empty after the first round. The probability that $J=j$ for any particular 
$j\in [m/8,7m/8]$ is
 \begin{eqnarray*}
  {m\choose j}\, \left(\frac{k}{r}\right)^{j}\, \left(\frac{r-k}{r}\right)^{m-{j}}
      &=& (1+o(1)) \, \sqrt{\frac{m}{2\pi {j}(m-{j})}}\, \left(\frac{mk}{{j}r}\right)^{j}\, \left(\frac{m(r-k)}{(m-{j})r}\right)^{m-{j}}\\
 &\le & (1+o(1)) \, \sqrt{\frac{m}{2\pi {j}(m-{j})}}\ =\ o(1),
   \end{eqnarray*}
 where we used Stirling's formula and the Arithmetic-Geometric Mean Inequality.
 On the other hand, we have whp that $m/8\le J\le 7m/8$ (by
Azuma's inequality and our assumption $k\in I$) and that the last $m$ balls
all go into different cells (since $m^2=o(r)$). Thus the probability of getting $R_i$ as the final ball distribution
is $o(1)$ uniformly in $i$, as desired.
 This finishes the second proof of Corollary~\ref{cr:grosu}.

\section*{Acknowledgements}

The author would like to thank Codrut Grosu for helpful comments.

\bibliography{oleg,general,ex,sets,graph}

\providecommand{\bysame}{\leavevmode\hbox to3em{\hrulefill}\thinspace}
\providecommand{\MR}{\relax\ifhmode\unskip\space\fi MR }
\providecommand{\MRhref}[2]{%
  \href{http://www.ams.org/mathscinet-getitem?mr=#1}{#2}
}
\providecommand{\href}[2]{#2}
\begin{thebibliography}{10}

\bibitem{alon+spencer:pm}
N.~Alon and J.~Spencer, \emph{The probabilistic method}, 3d ed., Wiley
  Interscience, 2008.

\bibitem{baber+talbot:11}
R.~Baber and J.~Talbot, \emph{Hypergraphs do jump}, Combin.\ Probab.\ Computing
  \textbf{20} (2011), 161--171.

\bibitem{brown+simonovits:84}
W.~G. Brown and M.~Simonovits, \emph{Digraph extremal problems, hypergraph
  extremal problems and the densities of graph structures}, Discrete Math.
  \textbf{48} (1984), 147--162.

\bibitem{erdos:64}
P.~Erd{\H{o}}s, \emph{On extremal problems of graphs and generalized graphs},
  Israel J.\ Math. \textbf{2} (1964), 183--190.

\bibitem{erdos+simonovits:66}
P.~Erd{\H{o}}s and M.~Simonovits, \emph{A limit theorem in graph theory},
  Stud.\ Sci.\ Math.\ Hungar. (1966), 51--57.

\bibitem{erdos+stone:46}
P.~Erd{\H o}s and A.~H. Stone, \emph{On the structure of linear graphs}, Bull.\
  Amer.\ Math.\ Soc. \textbf{52} (1946), 1087--1091.

\bibitem{frankl+rodl:84}
P.~Frankl and V.~R{\"o}dl, \emph{Hypergraphs do not jump}, Combinatorica
  \textbf{4} (1984), 149--159.

\bibitem{furedi:91}
Z.~F{\"u}redi, \emph{{Tur\'an} type problems}, Surveys in Combinatorics (A.~D.
  Keedwell, ed.), London Math.\ Soc.\ Lecture Notes Ser., vol. 166, Cambridge
  Univ.\ Press, 1991, pp.~253--300.

\bibitem{grosu:densities}
C.~Grosu, \emph{On the algebraic and topological structure of the set of
  {Tur\'an} densities}, E-print \texttt{arXiv:1403.4653v1}, 2014.

\bibitem{hardy+littlewood+polya:i}
G.~H. Hardy, J.~E. Littlewood, and G.~P{\'o}lya, \emph{Inequalities}, Cambridge
  Univ.\ Press, 1952, 2d ed.

\bibitem{katona+nemetz+simonovits:64}
G.~O.~H. Katona, T.~Nemetz, and M.~Simonovits, \emph{On a graph problem of
  {Tur\'an} {(In Hungarian)}}, Mat.\ Fiz.\ Lapok \textbf{15} (1964), 228--238.

\bibitem{keevash:11}
P.~Keevash, \emph{Hypergraph {Tur\'an} problem}, Surveys in Combinatorics
  (R.~Chapman, ed.), London Math.\ Soc.\ Lecture Notes Ser., vol. 392,
  Cambridge Univ.\ Press, 2011, pp.~83--140.

\bibitem{pikhurko:14}
O.~Pikhurko, \emph{Possible {Tur\'an} densities}, Israel J.\ Math. \textbf{201}
  (2014), 415--454.

\bibitem{sidorenko:95}
A.~Sidorenko, \emph{What we know and what we do not know about {Tur\'an}
  numbers}, Graphs Combin. \textbf{11} (1995), 179--199.

\bibitem{turan:41}
P.~{Tur\'an}, \emph{On an extremal problem in graph theory (in {Hungarian})},
  Mat.\ Fiz.\ Lapok \textbf{48} (1941), 436--452.

\end{thebibliography}
\end{document}